\documentclass[runningheads]{iwap}
\usepackage[margin=3.5cm]{geometry}
\usepackage{amsmath,amssymb}
\usepackage{url}

\newcommand{\lqn}[1]{\noalign{\noindent $\displaystyle{#1}$}}
\newcommand{\ind}{1\hspace{-.27em}\mbox{\rm l}}
\newcommand{\petitind}{1\hspace{-.2em}\mbox{\scriptsize\rm l}}

\begin{document}

\title*{A survey on the pseudo-process driven by the high-order heat-type
equation $\partial/\partial t=\pm\partial^N\!/\partial x^N$
concerning the hitting and sojourn times}
\toctitle{A survey on the pseudo-process driven by the high-order heat-type
equation $\partial/\partial t=\pm\partial^N\!/\partial x^N$}
\titlerunning{A pseudo-process driven by a high-order heat-type
equation}

\author{Aim\'e Lachal}
\index{Aim\'e Lachal}
\authorrunning{Aim\'e Lachal}
\institute{Institut Camille Jordan, UMR CNRS 5208 \\
Institut National des Sciences Appliqu\'ees de Lyon,\\
B\^at. L. de Vinci, 20 avenue Albert Einstein, 69621 Villeurbanne Cenex (France)\\
(e-mail: {aime.lachal@insa-lyon.fr})
}

\maketitle

\begin{abstract}

Fix an integer $N>2$ and let $X=(X(t))_{t\ge 0}$ be the pseudo-process driven
by the high-order heat-type equation $\partial/\partial t=\pm\partial^N\!/\partial x^N$.
The denomination ``pseudo-process'' means that $X$ is related to a signed
measure (which is not a probability measure) with total mass equal to 1.

In this survey, we present several explicit results and discuss some problems
concerning the pseudo-distributions of various functionals of the pseudo-process
$X$: the first or last overshooting times of a single barrier $\{a\}$ or a
double barrier $\{a,b\}$ by $X$; the sojourn times of $X$ in the intervals
$[a,+\infty)$ and $[a,b]$ up to a fixed time; the maximum or minimum of $X$
up to a fixed time.

\end{abstract}

\section{Introduction}

Consider the heat-type equation $\partial/\partial t=\kappa_{_{\! N}}
\partial^N\!/\partial x^N$ of order $N>2$ where $\kappa_{_{\! N}}=(-1)^{1+N/2}$
if $N$ is even and $\kappa_{_{\! N}}=\pm1$ if $N$ is odd. Let us introduce the
corresponding kernel $p(t;x)$ which is characterized by
$$
\int_{-\infty}^{+\infty} e^{iux} p(t;x)\, \mathrm{d}x
=\left\{\begin{array}{ll}
e^{-tu^N} &\mbox{if $N$ is even,} \\
e^{\kappa_{_{\! N}} t(-iu)^N} &\mbox{if $N$ is odd.}
\end{array}\right.
$$
This kernel defines a pseudo-process $(X(t))_{t\ge 0}$ driven by a signed
measure with total mass equal to 1 (which is not a probability measure)
according as the usual Markov rules: we set for $t>0$,
$0=t_0<t_1<\dots<t_m$ and $x=x_0,x_1,\dots,x_m,y\in\mathbb{R}$,
$$
\mathbb{P}_{\!x}\{X(t)\in \mathrm{d}y\}=p(t;x-y)\,\mathrm{d}y
$$
and
$$
\mathbb{P}_{\!x}\{X(t_1)\in \mathrm{d}x_1,\ldots,X(t_m)\in \mathrm{d}x_m\}=
\prod_{i=1}^{m}p(t_i-t_{i-1} ; x_{i-1}-x_i)\,\mathrm{d}x_i.
$$
Since we are dealing with a signed measure, it seems impossible
to extend the definition of the pseudo-process over all the positive
times. We can find in the literature two possible ad-hoc constructions:
one over the set of times of the form $kt/n$, $k,n\in\mathbb{N}$
(depending on a fixed time $t$, see~\cite{hoch} and~\cite{kry} for
pioneering works related to this construction),
the other one over the set of dyadic times $k/2^n$, $k,n\in\mathbb{N}$.
(which do not depend on any particular time, see~\cite{nish2} for this
last construction). For $N=2$, this is the most well-known Brownian motion
and for $N=4$, $(X(t))_{t\ge 0}$ is the so-called biharmonic pseudo-process.

For the pseudo-process $(X(t))_{t\ge 0}$ started at a point $x$,
we introduce:
\begin{itemize}
\item
the first overshooting times of a one-sided barrier $\{a\}$
(or, equivalently, the first hitting time of the half-line $[a,+\infty)$)
or a two-sided barrier $\{a,b\}$ (with the convention $\inf(\emptyset)=+\infty$):
\begin{align*}
\tau_a =\inf\{t\ge 0:X(t)\ge a\} \mbox{ for $x\le a$},\quad
\tau_{ab} =\inf\{t\ge 0:X(t)\notin[a,b]\}\mbox{ for $x\in[a,b]$};
\end{align*}

\item
the last overshooting times of such barriers before a fixed time $t$
(with the convention $\sup(\emptyset)=0$):
\begin{align*}
\sigma_a(t)=\sup\{s\in[0,t]:X(s)\ge a\},\quad \sigma_{ab}(t)
=\sup\{s\in[0,t]:X(t)\notin (a,b)\};
\end{align*}

\item
the sojourn times in the intervals $[a,+\infty)$ and $[a,b]$ up to a fixed
time $t$:
$$
T_a(t) =\mathrm{measure}\{s\in[0,t]:X(s)\ge a\},\quad
T_{ab}(t) =\mathrm{measure}\{s\in[0,t]:X(s)\in[a,b]\};
$$

\item
the maximum up to time~$t$:
$$
M(t)=\max_{0\le s\le t} X(s).
$$
\end{itemize}

In the foregoing rough definitions, the pseudo-distribution of the
quantity $T_a(t)$ for instance is to be understood as the limit of
$\mathbb{P}_{\!x}\{\frac1n \sum_{k=0}^{n-1}\ind_{[a,+\infty)}(X(kt/n))
\in \mathrm{d}s\}$ when $n\to\infty$.

We could introduce the alternative first hitting time of $(-\infty,a]$,
the alternative sojourn time in $(-\infty,a]$
and the up-to-date minimum $m(t)=\min_{0\le s\le t} X(s)$.
Actually, the pseudo-distributions of these three quantities are obviously
related to the pseudo-distributions of the foregoing ones.

We shall also consider the pseudo-process with a drift $(X^b(t))_{t\ge 0}$
defined by $X^b(t)=X(t)+bt$ where $b$ is a fixed real number. For this latter,
we introduce:
\begin{itemize}
\item
the first overshooting time of the threshold $a$:
$$
\tau_a^b=\inf\{t\ge 0:X^b(t)\ge a\}\mbox{ for $x\le a$}
$$
if the set $\{t\ge 0:X^b(t)\ge a\}$ is not empty, else we set
$\tau_a^b=+\infty$;

\item
the maximum functional up to time~$t$:
$$
M^b(t)=\max_{0\le s\le t} X^b(s).
$$
\end{itemize}

The aim of this survey is to provide a list of explicit results concerning
the pseudo-distributions of $(X(t),T_a(t))$, $(X(t),M(t))$,
$(\tau_a,X(\tau_a))$ and $\sigma_a(t)$, as well as those related to the
pseudo-process with a drift.
In particular, remarkable results hold for the pseudo-distributions of
$T_0(t)$ and $X(\tau_a)$. We also provide some methods for deriving those
of $T_{ab}(t)$, $\sigma_{ab}(t)$ and $(\tau_{ab},X(\tau_{ab}))$.

A way consists in using the Feynman-Kac functional
$$
\phi(t;x)=\mathbb{E}_x\!\!\left(e^{-\int_0^t f(X(s))\,\mathrm{d}s} g(X(t))\right)
\stackrel{\mbox{\scriptsize def}}{=} \lim_{n\to \infty}
\mathbb{E}_x\!\!\left(e^{-\frac tn\sum_{k=0}^{n-1} f(X(\frac{kt}{n}))}g(X(t))\right)
$$
which is a solution to the partial differential equation
$\frac{\partial \phi}{\partial t}(t;x)=\kappa_{_{\! N}}
\frac{\partial^N\!\phi}{\partial x^N}(t;x)-f(x)\phi(t;x)$ with
$\phi(0;x)=g(x)$. Its Laplace transform $\Phi(x)=\int_0^{+\infty} e^{-\lambda t}
\,\phi(t;x) \,\mathrm{d}t$ is a solution to the ordinary differential equation
$\kappa_{_{\! N}}\frac{\mathrm{d}^N\Phi}{\mathrm{d} x^N}(x)=(f(x)+\lambda)\Phi(x)-g(x).$
Another way consists in using Spitzer's identities which work actually when the
starting point is $0$ and $N$ is even. Indeed, their validity holds thanks
to the fact that the integral $\int_{-\infty}^{+\infty} |p(t;x)|\,\mathrm{d}x$
is finite, which is true only when $N$ is even. Additionally,
Spitzer's identities hinge on a symmetry property which is fulfilled only
when the starting point of the pseudo-process is $0$.
In the case $N=4$, see~\cite{nish3} for many
connections with fourth-order partial differential equations with various
boundary value conditions.

Let us introduce the $N^{\mbox{\scriptsize th}}$ roots of $\kappa_{_{\! N}}$:
$(\theta_\ell)_{1\le \ell\le N}$ and
$J=\{j\in\{1,\ldots,N\}:\Re e(\theta_j)>0\},$
$K=\{k\in\{1,\ldots,N\}:\Re e(\theta_k)<0\}$ that will be used
for solving the above differential equation. The notations $\# J$ and
$\#K$ stand for the cardinalities of the sets $J$ and $K$.
We have $\theta_\ell/\theta_m=e^{i(\ell-m)2\pi/N}$ for any $1\le\ell,m\le N$.

Set, for $j,j'\in J$ and $k,k'\in K$,
$$
A_j=\prod_{\ell\in J\setminus \{j\}} \frac{\theta_\ell}{\theta_\ell-\theta_j}
\quad\mbox{and}\quad
B_k=\prod_{\ell\in K\setminus \{k\}} \frac{\theta_\ell}{\theta_\ell-\theta_k}
,
$$
$$
C_{jj'k}=\prod_{j''\in J}(\theta_j\theta_{j'}-\theta_{j''}\theta_k)
\quad\mbox{and}\quad
D_{jkk'}=\prod_{k''\in K}(\theta_k\theta_{k'}-\theta_{k''}\theta_j).
$$
Let us also introduce the $(N-1)^{\mbox{\scriptsize th}}$ roots of
the complex number $i$:
$(\tilde{\theta}_l)_{1\le l\le N-1}$ and $\tilde{J}=\{j\in\{1,\ldots,N-1\}
:\Im(\tilde{\theta}_j)>0\}$. We shall need to introduce the roots
$(\omega_\ell^b(\lambda))_{1\le \ell\le N}$ of the polynomial $u^N+ibu+\lambda$
(where $\Re(\lambda)>0$). These last settings will be used for
the pseudo-process with a drift. Finally, set for any integer $\ell$ such
that $\ell\le N-1$
$$
I_{\ell}(t;\xi)=\frac{N i}{2 \pi} \left(e^{-i\frac{\ell}{N}\pi}
\int_0^{+\infty}\lambda^{N-\ell-1} e^{-t\lambda^N+e^{\frac{i\pi}{N}}
\xi\lambda} \,\mathrm{d}\lambda
-e^{i\frac{\ell}{N}\pi} \int_0^{+\infty}\lambda^{N-\ell-1}
e^{-t\lambda^N+e^{-\frac{i\pi}{N}}\xi\lambda} \,\mathrm{d}\lambda\right)\!.
$$
The functions $I_{\ell}$ satisfy $\int_0^{+\infty} e^{-\lambda t}
I_{\ell}(t;\xi)\,\mathrm{d}t=\lambda^{-\ell/N} e^{\!\!\sqrt[N\!]{\lambda}\,\xi}$
for $\lambda>0$ and $\Re(\xi)\le 0$.
They will be useful for expressing several distributions.

The results are presented by topic and in certain topics we have chosen to
exhibit them from the most particular to the most general thus following the
chronology. Moreover, it is not easy sometimes to deduce the particular cases
from the most general ones.
\section{Distributions related to $T_a(t)$}

See \cite{bho,bor,cl1,cl2,ho,kry,2003,2006,no,ors} for the chronology of
the results concerning the distributions related to $T_a(t)$ as well as for
the connections with the maximum and minimum functionals of $(X(t))_{t\ge0}$.

\subsection{Distribution of $T_a(t)$}

Set $\Phi(x)=\int_0^{+\infty} e^{-\lambda t} \,\mathbb{E}_x\!\!
\left(e^{-\mu T_a(t)}\right) \mathrm{d}t$ for $\lambda,\mu>0$ and $x\in\mathbb{R}$.
The quantity $\Phi(x)$ should be understood as
$$
\Phi(x)\stackrel{\mbox{\scriptsize def}}{=} \lim_{n\to \infty}
\int_0^{+\infty} e^{-\lambda t} \,\mathbb{E}_x\!\!\left(e^{-\mu\frac tn
\sum_{k=0}^{n-1} \petitind_{[a,+\infty)}(X(\frac{kt}{n}))}\right)\mathrm{d}t.
$$
Using the Feynman-Kac approach, it can be seen that the function $\Phi$ satisfies
the system
$$
\kappa_{_{\! N}}\frac{\mathrm{d}^N \Phi}{\mathrm{d} x^N}(x)
=\left\{\begin{array}{ll}
(\lambda+\mu)\,\Phi(x)-1  & \mbox{for } x\in(a,+\infty),
\\[1ex]
\lambda\, \Phi(x)-1 & \mbox{for } x\in(-\infty,a),
\end{array}\right.
$$
and
$$
\forall k\in\{0,1,\ldots,N-1\},\;\frac{\mathrm{d}^k \Phi}{\mathrm{d} x^k}(a^+)
=\frac{\mathrm{d}^k \Phi}{\mathrm{d} x^k}(a^-).
$$
This system can be explicitly solved by computing Vandermonde determinants.
In particular, for $x=a$, the following formula holds:
$$
\Phi(a)=\frac{1}{\!\sqrt[N\!\!]{\lambda^{\# K}(\lambda+\mu)^{\# J}}}
$$
and this two-parameters Laplace transform can be inverted (\cite{2003}).
The distribution of $T_a(t)$ under $\mathbb{P}_{\!a}$ is the same as that of
$T_0(t)$ under $\mathbb{P}_0$.
\begin{theorem}[Lachal, 2003]\label{th1}
The pseudo-distribution of $T_0(t)$ is a Beta law:
$$
\mathbb{P}_0\{T_0(t)\in \mathrm{d}s\}/\mathrm{d}s=\frac{1}{\pi}\,\sin\!
\Big(\frac{\# K}{N}\,\pi\Big)\,\frac{1}{\!\sqrt[N\!\!]{s^{\# K} (t-s)^{\# J}}}
\mbox{ for } s\in(0,t).
$$
\end{theorem}
\begin{example}
If $N$ is even, $T_0(t)$ obeys the famous Paul L\'evy's Arcsine law:
$$
\mathbb{P}_0\{T_0(t)\in \mathrm{d}s\}/\mathrm{d}s=\frac{1}{\pi\sqrt{s(t-s)}}.
$$
In the history of pseudo-processes, this result was discovered by Krylov (\cite{kry})
when $N$ is even. For $N=3$, Orsingher obtained (\cite{ors})
$$
\mathbb{P}_0\{T_0(t)\in \mathrm{d}s\}/\mathrm{d}s
=\begin{cases}
\displaystyle\frac{\sqrt3}{2\pi\sqrt[3]{s^2(t-s)}} & \mbox{when $\kappa_{_{\! N}}=+1$,}
\\[2ex]
\displaystyle\frac{\sqrt3}{2\pi\sqrt[3]{s(t-s)^2}} & \mbox{when $\kappa_{_{\! N}}=-1$.}
\end{cases}
$$
\end{example}

Using a similar method, the following simple results can be obtained (\cite{2003}).
\begin{theorem}[Lachal, 2003]\label{th2}
The pseudo-distribution of $T_0(t)$ conditioned on $X(t)=0$ is the uniform law
on $(0,t)$: for $s\in(0,t)$,
$$
\mathbb{P}_0\{T_0(t)\in \mathrm{d}s|X(t)=0\}/\mathrm{d}s=\frac{1}{t}.
$$
The pseudo-distributions of $T_0(t)$ conditioned on $X(t)>0$ and $X(t)<0$
are Beta laws: for $s\in(0,t)$,
\begin{align*}
\mathbb{P}_0\{T_0(t)\in \mathrm{d}s|X(t)>0\}/\mathrm{d}s
&
=\frac{N\sin(\frac{\#K}{N}\,\pi)}{(\#K\pi)t}
\left(\frac{s}{t-s}\right)^{\!\!\frac{\#J}{N}},
\\
\mathbb{P}_0\{T_0(t)\in \mathrm{d}s|X(t)<0\}/\mathrm{d}s
&
=\frac{N\sin(\frac{\#J}{N}\,\pi)}{(\#J\pi)t}
\left(\frac{t-s}{s}\right)^{\!\!\frac{\#K}{N}}.
\end{align*}
\end{theorem}
\begin{example}
If $N$ is even, for $s\in(0,t)$,
\begin{align*}
\mathbb{P}_0\{T_0(t)\in \mathrm{d}s|X(t)>0\}/\mathrm{d}s
=\frac{2}{\pi t} \sqrt{\frac{s}{t-s}},\quad
\mathbb{P}_0\{T_0(t)\in \mathrm{d}s|X(t)<0\}/\mathrm{d}s
=\frac{2}{\pi t} \sqrt{\frac{t-s}{s}}.
\end{align*}
\end{example}

The results of Theorems~\ref{th1} and~\ref{th2} were found by Hochberg,
Nikitin and Orsingher (\cite{ho,no,ors}) in the cases $N=3,4,5,7$ and
conjectured in the general case.

\subsection{Distribution of $(X(t),T_a(t))$}

\subsubsection{\textsl{Case $x=a$.}}
Set $\Phi=\int_0^{\infty} e^{-\lambda t} \mathbb{E}_0(e^{i \mu X(t)-\nu T_0(t)})
\,\mathrm{d}t$ for $\lambda,\nu>0$ and $\mu\in\mathbb{R}$. The quantity
$\Phi$ can be understood as
\begin{align*}
\Phi
&\stackrel{\mbox{\scriptsize def}}{=}
\lim_{n\to \infty} \sum_{k=0}^{\infty} \int_{k/2^n}^{(k+1)/2^n}
e^{-\lambda t} \,\mathbb{E}_0\!\!\left(e^{i\mu X(k/2^n)-\frac{\nu}{2^n}
\sum_{j=1}^{k}\petitind_{[0,+\infty)}(X(j/2^n))}\right)\mathrm{d}t
\\
&
=\lim_{n\to\infty} \frac{1-e^{-\lambda/2^n}}{\lambda} \sum_{k=0}^{\infty}
e^{-\lambda k/2^n} \mathbb{E}_0\!\!\left(e^{i\mu X(k/2^n)-\frac{\nu}{2^n}
\sum_{j=1}^{k}\petitind_{[0,+\infty)}(X(j/2^n))}\right)\!.
\end{align*}
A Spitzer's identity yields the following relationship which holds for $|z|,|\zeta|<1$:
\begin{align*}
\lqn{\sum_{k=0}^{\infty} \mathbb{E}_0\!\Big(e^{i \mu X(k/2^n)}
\zeta^{\sum_{j=1}^{k}\petitind_{[0,+\infty)}(X(j/2^n))}\Big) z^k}
\\[-4ex]
&=\frac{1}{1-z}\,\exp\!\left(\,\sum_{k=1}^{\infty} \frac{1}{k} \,\mathbb{E}_0\!
\Big(e^{i\mu X(k/2^n)}\zeta^{k\petitind_{[0,+\infty)}(X(k/2^n))}-1\Big)z^k \right)\!.
\end{align*}
With this identity at hand, it can be seen that
$$
\Phi=\frac{1}{\displaystyle\mathop{\textstyle{\prod}}_{j \in J}
(\!\sqrt[N\!]{\lambda+\nu}-i\mu \theta_j )\displaystyle\mathop{
\textstyle{\prod}}_{k \in K} (\!\sqrt[N\!]{\lambda}-i \mu \theta_k )}.
$$
This three-parameters Laplace-Fourier transform can be inverted (\cite{cl1}).
\begin{theorem}[Cammarota \& Lachal, 2010]
The pseudo-distribution of the vector $(X(t),T_0(t))$ is given, for $s\in(0,t)$
and
$y \le 0$, by
\begin{align*}
\lqn{\mathbb{P}_0\{X(t) \in \mathrm{d}y, T_0(t) \in \mathrm{d}s\}/
\mathrm{d}y \,\mathrm{d}s}
&
=-\frac{Ni}{2\pi}\sum_{m=0}^{\# K} \alpha_{-m} \,s^{\frac{m-\#K}{N}}
\int_0^{\infty} \xi^{m+\# J} e^{-(t-s)\xi^N}
\mathcal{K}_m(y\xi) \,E_{1,\frac{m+\# J}{N}}(-s \xi^N)\,\mathrm{d}\xi
\end{align*}
and, for $s\in(0,t)$ and $y \ge 0$, by
\begin{align*}
\lqn{\mathbb{P}_0\{X(t)\in \mathrm{d}y, T_0(t)\in \mathrm{d}s\}/
\mathrm{d}y \,\mathrm{d}s}
&
=\frac{Ni}{2 \pi} \sum_{m=0}^{\# J} \beta_{-m}\, (t-s)^{\frac{m-\#J}{N}}
\int_0^{\infty} \xi^{m+\# K} \,e^{-s \xi^N} \,\mathcal{J}_m(y\xi)\,
E_{1,\frac{m+\# K}{N}} \!\left(-(t-s)\xi^N\right) \mathrm{d}\xi
\end{align*}
where $\alpha_m=\sum_{j\in J} A_j\theta_j^m$, $\beta_m=\sum_{k\in K} B_k\theta_k^m$
for any integer $m$ (in particular $\alpha_0=\beta_0=1$, $\alpha_{-\#K}=1,
\beta_{-\#J}=(-1)^{\#J}$),
\begin{align*}
\mathcal{J}_m(z)
&=e^{-i \frac{\# J-m-1}{N}\pi} \sum_{j \in J} A_j \theta_j^{m+1}
e^{-\theta_j e^{i\frac{\pi}{N}}z}-e^{i \frac{\# J-m-1}{N}\pi}
\sum_{j \in J} A_j \theta_j^{m+1} e^{-\theta_j e^{-i\frac{\pi}{N}}z},
\\
\mathcal{K}_m(z)
&=e^{-i \frac{\# K-m-1}{N}\pi} \sum_{k \in K} B_k\theta_k^{m+1}
e^{-\theta_k e^{i\frac{\pi}{N}}z}-e^{i \frac{\# K -m-1}{N}\pi}
\sum_{k \in K} B_k\theta_k^{m+1}e^{-\theta_k e^{-i\frac{\pi}{N}}z},
\end{align*}
and $E_{a,b}$ is the Mittag-Leffler function $E_{a,b}(z)=\sum_{n=0}^{\infty}
\frac{z^n}{\Gamma(an+b)}$.
\end{theorem}
%
\begin{remark}
By choosing $y=0$ in the pseudo-distribution of $(X(t),T_0(t))$ in the
foregoing theorem and next dividing the result by $\mathbb{P}_{\!x}\{X(t)=0\}$,
we could retrieve the pseudo-distribution of the corresponding sojourn time
of the ``pseudo-bridge'' $(X(s)|X(0)=X(t)=0)_{0\le s\le t}$ displayed in
Theorem~\ref{th2}.
\end{remark}
%
\subsubsection{\textsl{Case $x\ne a$.}}

Set $\Phi(x,y)=\int_0^{+\infty} e^{-\lambda t}
\big(\mathbb{E}_x \big(e^{-\mu T_a(t)},X(t)\in \mathrm{d}y\big)/\mathrm{d}y\big) \,\mathrm{d}t$
for $\lambda,\mu>0$ and $x,y\in\mathbb{R}$. It can be seen that $\Phi$ solves the
differential equation
$$
\kappa_{_{\! N}}\frac{\partial^N \Phi}{\partial x^N}(x,y)
=\left\{\begin{array}{ll}
(\lambda+\mu)\,\Phi(x,y)-\delta_y(x)  & \mbox{for } x\in(a,+\infty),
\\[1ex]
\lambda\, \Phi(x,y)-\delta_y(x) & \mbox{for } x\in(-\infty,a),
\end{array}\right.
$$
with regularity conditions
$$
\left\{\begin{array}{l}
\displaystyle\forall k\in\{0,1,\ldots,N-1\},\;
\frac{\partial^k \Phi}{\partial x^k}(a^+,y)=\frac{\partial^k \Phi}{\partial x^k}(a^-,y),
\\[2ex]
\displaystyle\forall k\in\{0,1,\ldots,N-2\},\;
\frac{\partial^k \Phi}{\partial x^k}(y^+,y)=\frac{\partial^k \Phi}{\partial x^k}(y^-,y)
\mbox{ and }\displaystyle\frac{\partial^{N-1}\Phi}{\partial x^{N-1}}(y^+,y)-
\frac{\partial^{N-1}\Phi}{\partial x^{N-1}}(y^-,y)=\kappa_{_{\! N}}.
\end{array}\right.
$$
This system can be explicitly solved by computing Vandermonde determinants
and the inversion of the two-parameters Laplace transform can be performed
(\cite{cl2}).

\begin{theorem}[Cammarota \& Lachal, 2010]
Set
\begin{eqnarray*}
f(t;\theta)=\frac{\sin(\frac \pi N)}{\pi N}\frac{t^{\frac 1 N-1}(\theta^2
-t^\frac 2 N)}{(t^{\frac 2N}-2\theta\cos(\frac\pi N)t^{\frac 1N}+\theta^2)^2}.
\end{eqnarray*}
1. Assume that $y\ge 0$. For $s\in(0,t)$, if $x\in(-\infty,0]$,
\begin{align*}
\lqn{\mathbb{P}_{\!x}\{X(t)\in\mathrm{d}y, T_0(t)\in\mathrm{d}s\}/
\mathrm{d}y\,\mathrm{d}s}
&\\[-4ex]
&
=\kappa_{_{\! N}} \ind_{\{\# J=\#K+1\}} \Bigg(\sum_{j \in J}\theta_j
A_j I_{\# K}(t-s;-\theta_j x)\Bigg)
\!\Bigg(\sum_{j\in J} \theta_j^{\# J} A_j I_{\#K}(s;\theta_j y) \Bigg)
\\
&
\hphantom{=\;}+\kappa_{_{\! N}} \sum_{j,j'\in J,\,k\in K}
\frac{A_j\theta_{j'}A_{j'}\theta_k^{\#K-\#J}
B_k C_{j j' k}}{\theta_j^{\#K-1}(\theta_{j'}-\theta_k )}
\int_0^s \sigma^{\frac 1 N-1} I_{\#J-1}(s-\sigma;\theta_{j'} y)
\,\mathrm{d}\sigma
\\
&
\hphantom{=\;}\times \int_0^{t-s} I_{\#K}(\tau;-\theta_j x)
f\Big(t-s-\tau;\frac{\theta_j}{\theta_k}\!\sqrt[N\!] \sigma\,\Big) \,\mathrm{d}\tau
\end{align*}
and if $x\in[0,\infty)$,
\begin{align*}
\lqn{\mathbb{P}_{\!x}\{X(t)\in\mathrm{d}y, T_0(t)\in\mathrm{d}s\}/
\mathrm{d}y\,\mathrm{d}s}&
\\[-4ex]
&
=\kappa_{_{\! N}} \sum_{j,j'\in J,\,k \in K} \frac{\theta_j A_j A_{j'}
B_k C_{jj'k}}{\theta_k^{\#J-1}(\theta_{j}-\theta_k)}
\int_0^s \sigma^{\frac 1N-1}I_{N-2} (s-\sigma;\theta_{j}y-\theta_k x)
\,\mathrm{d}\sigma
\\
&
\hphantom{=\;}\times
\int_0^{t-s} \frac{(t-s-\tau)^{\frac 1 N-1}}{\Gamma(\frac 1 N)}
f\Big(\tau; \frac {\theta_{j'}}{ \theta_k}\!\sqrt[N\!] \sigma\,\Big)
\,\mathrm{d}\tau.
\end{align*}
For $s=t$, if $x\in(-\infty,0]$,
$\mathbb{P}_{\!x}\{X(t)\in\mathrm{d}y,T_0(t)=t\}/ \mathrm{d}y=0$
and if $x\in[0,+\infty)$, there is an atom given by
$$
\mathbb{P}_{\!x}\{X(t)\in\mathrm{d}y,T_0(t)=t\}/ \mathrm{d}y
=p(t;x-y)+ \sum_{j\in J,\,k\in K}\frac{\theta_j A_j
\theta_k B_k}{\theta_j-\theta_k} \,I_{N-1}(t;\theta_j y-\theta_k x).
$$
2. Assume that $y\le 0$. For $s\in(0,t)$, if $x\in(-\infty,0]$,
\begin{align*}
\lqn{\mathbb{P}_{\!x}\{X(t)\in\mathrm{d}y, T_0(t)\in\mathrm{d}s\}/
\mathrm{d}y\,\mathrm{d}s}&
\\[-4ex]
&
=-\sum_{j\in J,\,k,k'\in K} \frac{A_j\theta_k B_kB_{k'} D_{j k k'}}{\theta_j^{\#K-2}
\theta_{k'}(\theta_k-\theta_j)} \,s^{\frac 1N-1} \int_0^{t-s}
I_{N-1}(\tau;\theta_k y-\theta_j x)
f\Big(t-s-\tau;\frac{\theta_{j}}{\theta_{k'}} \!\sqrt[N\!]s\,\Big)\,\mathrm{d}\tau
\end{align*}
and if $x\in [0,\infty)$,
\begin{align*}
\lqn{\mathbb{P}_{\!x}\{X(t)\in\mathrm{d}y, T_0(t)\in\mathrm{d}s\}/ \mathrm{d}y\,\mathrm{d}s}&
\\[-4ex]
&
=\ind_{\{\# K=\#J+1\}} \Bigg(\sum_{k\in K} \theta_k B_k I_{\# J}(s;-\theta_k x)\Bigg)
\!\Bigg(\sum_{k\in K} \theta_k^{\# K} B_k I_{\#J}(t-s;\theta_k y)\Bigg)
\\
&
\hphantom{=\;} +\sum_{j\in J,\,k,k'\in K} \frac{\theta_j^{\#J-\#K+1}A_j B_k
\theta_{k'}B_{k'} D_{jkk'}}{\theta_k^{\#J} (\theta_j-\theta_{k'})}
\int_0^s \sigma^{\frac 1N-1} I_{\#J-1}(s -\sigma;-\theta_k x) \,\mathrm{d}\sigma
\\
&
\hphantom{=\;}\times \int_0^{t-s} I_{\#K}(\tau;\theta_{k'} y)
f\Big(t-s-\tau;\frac{\theta_j}{\theta_k}\!\sqrt[N\!] \sigma\,\Big) \,\mathrm{d}\tau.
\end{align*}
For $s=0$, if $x\in [0,\infty)$, $\mathbb{P}_{\!x}\{X(t)\in\mathrm{d}y,T_0(t)=0\}/\mathrm{d}y=0$
and if $x\in(-\infty,0]$, there is an atom given by
$$
\mathbb{P}_{\!x}\{X(t)\in\mathrm{d}y,T_0(t)=0\}/ \mathrm{d}y
=p(t;x-y)+\sum_{j\in J,\,k\in K} \frac{\theta_j A_j
\theta_k B_k}{\theta_j-\theta_k} \,I_{N-1}(t;\theta_k y-\theta_j x).
$$
\end{theorem}
The following relationship between $T_0(t)$, $M(t)$ and $\tau_0$ holds:
for $x,y\le 0$,
$$
\mathbb{P}_{\!x}\{X(t)\in\mathrm{d}y,T_0(t)=0\}/ \mathrm{d}y
=\mathbb{P}_{\!x}\{X(t)\in\mathrm{d}y,M(t)< 0\}/ \mathrm{d}y
=\mathbb{P}_{\!x}\{X(t)\in\mathrm{d}y,\tau_0>t\}/ \mathrm{d}y.
$$
%
\begin{remark}
By integrating the joint pseudo-distribution of $(X(t),T_0(t))$ with respect to $y$
in the foregoing theorem, we could derive the pseudo-distribution of $T_0(t)$.
Actually, the result does not simplify so much.
\end{remark}

\section{Distribution related to $M(t)$}

See~\cite{bho,bor,hoch,2003,2006,2007,2008,nish2} for references
related to the various distributions related to $M(t)$.

\subsection{Distribution of $M(t)$}

The variables $T_a(t)$ and $M(t)$ are related together according to
$$
\mathbb{P}_{\!x}\{M(t)\le a\}=\mathbb{P}_{\!x}\{T_a(t)=0\}=\lim_{\mu\to+\infty}
\mathbb{E}_x\!\big(e^{-\mu T_a(t)}\big).
$$
The quantity $\mathbb{P}_{\!x}\{M(t)\le a\}$ should be understood as
$$
\mathbb{P}_{\!x}\{M(t)\le a\}\stackrel{\mbox{\scriptsize def}}{=}
\lim_{n\to \infty} \mathbb{P}_{\!x}\bigg\{\max_{0\le k\le n}
X\!\bigg(\frac{kt}{n}\bigg)\!\le a\bigg\}.
$$
Thanks to this connection, it is possible to deduce that, for $x\le a$
and $\lambda>0$,
$$
\int_0^{+\infty}  e^{-\lambda t} \,\mathbb{P}_{\!x} \{M(t)\le a\} \,\mathrm{d}t
=\frac{1}{\lambda}\Bigg(1-\sum_{j\in J} A_j \, e^{\theta_k\!\!\sqrt[N\!]{\lambda}
\,(x-a)}\Bigg)\!.
$$
The foregoing Laplace transform can be inverted (\cite{2003}).
\begin{theorem}[Lachal, 2003]
The pseudo-distribution of $M(t)$ is given by
\begin{align*}
\mathbb{P}_{\!x}\{M(t)\ge a\}
&
=\sum_{m=0}^{\#J-1} a_m \int_0^t \frac{\partial^m p}{\partial x^m}(s;x-a)
\,\frac{\mathrm{d}s}{(t-s)^{1-(m+1)/N}}
\\
&
=2\,\mathbb{P}_{\!x}\{X(t)\ge a\}-\sum_{m=0}^{\#J-1} b_m \int_0^t \frac{\partial^m p}{\partial x^m}(s;x-a)
\,\frac{\mathrm{d}s}{(t-s)^{1-(m+1)/N}}
\end{align*}
where
$$
a_m=\frac{(-1)^m N}{\Gamma(\frac{m+1}{N})}\sum_{j\in J} A_j^2\,\sigma_{j,\#J-1-m}
\quad\mbox{ and }\quad b_m=\frac{2(-1)^m}{\Gamma(\frac{m+1}{N})}\sum_{j\in J} A_j\,\sigma_{j,\#J-1-m}-a_m.
$$
The coefficients $\sigma_{j,p}$, $j\in J$, $0\le p\le\#J-1$, are given by
$\sigma_{j,0}=1$ and for $1\le p\le \#J-1$,
$$
\sigma_{j,p}=\sum_{\ell_1,\ldots,\ell_p\in J\setminus\{j\}\atop
\ell_1<\cdots<\ell_p}\theta_{\ell_1}\cdots\theta_{\ell_p}.
$$
\end{theorem}
%
\begin{remark}
From the second displayed expression of $\mathbb{P}_{\!x}\{M(t)\ge a\}$,
we can see that the famous reflection principle for Brownian motion does
not hold any longer for pseudo-processes related to an order $N>2$.
\end{remark}
%
\begin{example}
For $N=3$, Orsingher (\cite{ors}) derived the historical result
\begin{align*}
\lqn{\mathbb{P}_{\!x}\{M(t)\ge a\}}
=\begin{cases}
\displaystyle\frac{3}{\Gamma(\frac13)}\int_0^t p(s;x-a)
\,\frac{\mathrm{d}s}{(t-s)^{\frac23}} & \mbox{when $\kappa_{_{\! N}}=+1$,}
\\[2ex]
\displaystyle\frac{2}{\Gamma(\frac13)}\int_0^t p(s;x-a)
\,\frac{\mathrm{d}s}{(t-s)^{\frac23}}-\,\frac{1}{\Gamma(\frac23)}
\int_0^t \frac{\partial p}{\partial x}(s;x-a)\,\frac{\mathrm{d}s}{(t-s)^{\frac13}}
& \mbox{when $\kappa_{_{\! N}}=-1$.}
\end{cases}
\end{align*}
For $N=4$, Hochberg (\cite{hoch}) derived the historical result
$$
\int_0^{+\infty} e^{-\lambda t}\,(\mathbb{P}_{\!x}\{M(t)\in \mathrm{d}a\}/
\mathrm{d}a)\,\mathrm{d}t
=-\frac{\sqrt 2}{\lambda^{3/4}} \,e^{\sqrt[4]{\lambda} \,(x-a)/\sqrt 2}
\sin \!\bigg(\frac{\sqrt[4]{\lambda} \,(x-a)}{\sqrt 2}\bigg)
$$
which was subsequently completed by Beghin, Orsingher and Ragozina (\cite{bor}):
$$
\mathbb{P}_{\!x}\{M(t)\ge a\}=\frac{2\sqrt2}{\Gamma(\frac14)}
\int_0^t p(s;x-a)\,\frac{\mathrm{d}s}{(t-s)^{\frac34}}.
$$
\end{example}
\subsection{Distribution of $(X(t),M(t))$}

Set $\Phi(x)=\int_0^{+\infty} e^{-\lambda t} \,\mathbb{E}_x
\big(e^{i\mu X(t)-\nu M(t)}\big)\,\mathrm{d}t$ for $\lambda,\nu>0$,
$\mu\in\mathbb{R}$ and $x\in\mathbb{R}$.
The quantity $\Phi(x)$ can be understood as
\begin{align*}
\Phi(x)
&\stackrel{\mbox{\scriptsize def}}{=}
\lim_{n\to \infty} \sum_{k=0}^{\infty} \int_{k/2^n}^{(k+1)/2^n}
e^{-\lambda t} \,\mathbb{E}_x\!\!\left(e^{i\mu X(k/2^n)-
\nu\max_{0\le j\le k}X(j/2^n)}\right)\mathrm{d}t
\\
&
=\lim_{n\to\infty} e^{(i\mu-\nu)x}\frac{1-e^{-\lambda/2^n}}{\lambda}
\sum_{k=0}^{\infty}e^{-\lambda k/2^n} \mathbb{E}_0\!\!\left(e^{i\mu X(k/2^n)
-\nu\max_{0\le j\le k}X(j/2^n)}\right)\!.
\end{align*}
Another Spitzer's identity yields the following relationship which holds for $|z|<1$:
$$
\sum_{k=0}^{\infty} \mathbb{E}_0\Big(e^{i\mu X(k/2^n)-
\nu \max_{0\le j\le k}X(j/2^n)}\Big)z^k
=\frac{1}{1-z}\,\exp\!\left(\,\sum_{k=1}^{\infty} \frac{1}{k} \,\mathbb{E}_0\!
\Big(e^{i\mu X(k/2^n)-\nu X(k/2^n)^+}\!-1\Big)z^k \right)\!.
$$
The Laplace-Fourier transform of the vector $(X(t),M(t))$ ensues:
$$
\Phi(x)=\frac{e^{(i\mu-\nu)x}}{\displaystyle\mathop{\textstyle{\prod}}_{j\in J}
(\!\sqrt[N\!]{\lambda}-(i\mu-\nu)\theta_j)
\displaystyle\mathop{\textstyle{\prod}}_{k\in K} (\!\sqrt[N\!]{\lambda}-i\mu\theta_k)}.
$$
This three-parameters transform can be progressively inverted (\cite{2007}).
For $z\ge x\vee y$,
$$
\int_0^{+\infty} e^{-\lambda t}[\mathbb{P}_{\!x}\{X(t)\in \mathrm{d}y,M(t)
\in \mathrm{d}z\}/\mathrm{d}y\,\mathrm{d}z]\,\mathrm{d}t
=\frac{1}{\lambda}\,\chi_{_{\scriptstyle J}}(\lambda;x-z)\,
\chi_{_{\scriptstyle K}}(\lambda;z-y),
$$
with
$$
\chi_{_{\scriptstyle J}}(\lambda;\xi)=\!\sqrt[N\!]{\lambda}
\,\sum_{j\in J} \theta_jA_j\,e^{\theta_j\!\!\sqrt[N\!]{\lambda}\,\xi},
\quad
\chi_{_{\scriptstyle K}}(\lambda;\xi)=-\!\sqrt[N\!]{\lambda}
\,\sum_{k\in K} \theta_kB_k\,e^{\theta_k\!\!\sqrt[N\!]{\lambda}\,\xi}.
$$
This last Laplace transform can also be inverted (\cite{2007}).
\begin{theorem}[Lachal, 2006]
The joint pseudo-distribution of $(X(t),M(t))$ admits the representation below.
For $z\ge x\vee y$,
\begin{align*}
\mathbb{P}_{\!x}\{X(t)\le y\le z\le M(t)\}
&=\sum_{k\in K\atop 0\le m\le\#J-1}
a_{km}\int_0^t\!\!\int_0^s\frac{\partial^m p}{\partial x^m}(\sigma;x-z)
\,\frac{I_{0}(s-\sigma;\theta_k(z-y))}{(t-s)^{1-(m+1)/N}}\,\mathrm{d}s\,
\mathrm{d}\sigma
\\
&=\sum_{j\in J\atop 0\le m\le\#K-1}
b_{jm}\int_0^t\!\!\int_0^s\frac{\partial^m p}{\partial x^m}(\sigma;z-y)
\,\frac{I_{0}(s-\sigma;\theta_j(x-z))}{(t-s)^{1-(m+1)/N}}\,\mathrm{d}s\,
\mathrm{d}\sigma
\end{align*}
where
$$
a_{km}=\frac{(-1)^m NB_k}{\Gamma(\frac{m+1}{N})}\sum_{j\in J}
\frac{\theta_jA_j^2\sigma_{j,\#J-1-m}}{\theta_j-\theta_k}
\quad\mbox{and}\quad b_{jm}=\frac{(-1)^m N\theta_jA_j}{\Gamma(\frac{m+1}{N})}\sum_{k\in K}
\frac{B_k^2\sigma_{k,\#K-1-m}}{\theta_k-\theta_j}.
$$
\end{theorem}

\begin{example}
For $N=3$, Beghin, Orsingher and Ragozina (\cite{bor}) derived the result
\begin{align*}
\lqn{\mathbb{P}_{\!x}\{X(t)\le y\le z\le M(t)\}}
=\begin{cases}
\displaystyle\frac{1}{\Gamma(1/3)}\int_0^t\!\!\int_0^s p(\sigma;x-z)\,q(s-\sigma;z-y)\,
\frac{\mathrm{d}s\,\mathrm{d}\sigma}{(t-s)^{2/3}} & \mbox{when $\kappa_{_{\! N}}=+1$,}
\\[2ex]
\displaystyle\frac{1}{\Gamma(1/3)}\int_0^t\!\!\int_0^s
p(\sigma;z-y)\,q(s-\sigma;x-z)\, \frac{\mathrm{d}s\,\mathrm{d}\sigma}{(t-s)^{2/3}}
& \mbox{when $\kappa_{_{\! N}}=-1$,}
\end{cases}
\end{align*}
with
$$
p(t;\xi)=\frac{1}{\pi} \int_0^{+\infty} \cos(\xi\lambda-t\lambda^3)\,\mathrm{d}\lambda
$$
and
\begin{align*}
q(t;\xi)&=\begin{cases}
\displaystyle
\frac{\xi}{\pi t} \int_0^{+\infty} e^{-t\lambda^3+\frac12\,\xi\lambda}
\sin\!\bigg(\frac{\sqrt3}{2}\,\xi\lambda+\frac{\pi}{3}\bigg)\,\mathrm{d}\lambda
& \mbox{when $\kappa_{_{\! N}}=+1$,}
\\[3ex]
\displaystyle
\frac{\xi}{\pi t} \bigg[\sqrt3 \int_0^{+\infty} e^{-t\lambda^3+\xi\lambda}\,\mathrm{d}\lambda
+\int_0^{+\infty} e^{-t\lambda^3-\frac12\,\xi\lambda}
\sin\!\bigg(\frac{\sqrt3}{2}\,\xi\lambda+\frac{\pi}{3}\bigg)\,\mathrm{d}\lambda\bigg]
& \mbox{when $\kappa_{_{\! N}}=-1$.}
\end{cases}
\end{align*}
For $N=4$, they derived
\begin{align*}
\mathbb{P}_{\!x}\{X(t)\le y\le z\le M(t)\}
&=\int_0^t\!\!\int_0^s p(\sigma;x-z)\,q_1(s-\sigma;z-y)\, \frac{\mathrm{d}s\,\mathrm{d}\sigma}{(t-s)^{3/4}}
\\
&
\hphantom{=\,}
+\int_0^t\!\!\int_0^s \frac{\partial p}{\partial x}(\sigma;x-z)\,q_2(s-\sigma;z-y)\,
\frac{\mathrm{d}s\,\mathrm{d}\sigma}{\sqrt{t-s}}
\end{align*}
with
\begin{align*}
q_1(t;\xi)
&=\frac{\xi}{\pi\sqrt2\,\Gamma(1/4)\,t} \int_0^{+\infty} e^{-t\lambda^4}\,\cos(\xi\lambda)\,\mathrm{d}\lambda,
\\[1ex]
q_2(t;\xi)
&=\frac{\xi}{2\pi^2\,t} \int_0^{+\infty} e^{-t\lambda^4}
\left[\cos(\xi\lambda)+\sin(\xi\lambda)-e^{-\xi\lambda}\right]\mathrm{d}\lambda.
\end{align*}
\end{example}

\subsection{Distribution of $(X^b(t),M^b(t))$}

The Laplace-Fourier transform of the vector $(X^b(t),M^b(t))$ is given,
for $\lambda,\nu>0$ and $\mu\in\mathbb{R}$, by
$$
\mathbb{E}_x\!\!\left(\int_0^{+\infty} e^{-\lambda t+i\mu X^b(t)-\nu M^b(t)}
\,\mathrm{d}t\right)=\frac{e^{(i\mu-\nu)x}}{\displaystyle\mathop{
\textstyle{\prod}}_{j\in J}(\omega_j^b(\lambda)+\mu+i\nu)
\displaystyle\mathop{\textstyle{\prod}}_{k\in K} (\omega_k^b(\lambda)+\mu)}.
$$
This three-parameters transform can be partially inverted for giving the
following result (\cite{2008}).
\begin{theorem}[Lachal, 2008]
The Laplace transform with respect to time $t$ of the joint pseudo-distribution of
$(X^b(t),M^b(t))$ is given, for $\lambda>0$ and $z\ge x\vee y$, by
$$
\int_0^{+\infty} e^{-\lambda t}[\mathbb{P}_{\!x}\{X^b(t)\in
\mathrm{d}y,M^b(t)\in \mathrm{d}z\}/\mathrm{d}y\,\mathrm{d}z]\,\mathrm{d}t
=\chi_{_{\scriptstyle J}}^b(\lambda;x-z)\,\chi_{_{\scriptstyle K}}^b(\lambda;z-y)
$$
where
$$
\chi_{_{\scriptstyle J}}^b(\lambda;\xi)= \sum_{j\in J}
\frac{e^{-i\omega_j^b(\lambda)\xi}}{\displaystyle\mathop{
\textstyle{\prod}}_{\ell\in J\setminus\{j\}}
(\omega_\ell^b(\lambda)-\omega_j^b(\lambda))}
\quad\mbox{and}\quad \chi_{_{\scriptstyle K}}^b(\lambda;\xi)
=\sum_{k\in K} \frac{e^{-i\omega_k^b(\lambda)\xi}}%
{\displaystyle\mathop{\textstyle{\prod}}_{l\in K\setminus\{k\}}
(\omega_\ell^b(\lambda)-\omega_k^b(\lambda))}.
$$
\end{theorem}

\subsection{Distribution of $\sigma_a(t)$}

The variables $\sigma_a(t)$ and $M(t)$ are related together according as
$$
\mathbb{P}_{\!x}\{\sigma_a(t)\le s\}=\mathbb{P}_{\!x}\Big\{
\max_{s\le u\le t} X(u)\le a\Big\}=\int_{-\infty}^a p(s;x-y)\,
\mathbb{P}_{\!y}\{ M(t-s)\le a\} \,\mathrm{d}y
$$
The pseudo-distribution of $\sigma_a(t)$ under $\mathbb{P}_{\!a}$ is the same
as that of $\sigma_0(t)$ under $\mathbb{P}_0$.
\begin{theorem}[Lachal, 2003]
The iterated Laplace transform of $\sigma_0(t)$ under $\mathbb{P}_0$ is given,
for $\lambda,\mu>0$, by
$$
\int_0^{+\infty} e^{-\lambda t}\,\mathbb{E}_0\big(e^{-\mu \sigma_0(t)}\big)
\,\mathrm{d}t=\frac{1}{\lambda+\mu}\Bigg[1-\frac1N \sum_{j\in J}\prod_{j'\in J}
\left(1-\frac{\theta_j}{\theta_{j'}}\!\sqrt[N\!]{\frac{\lambda+\mu}{\lambda}}
\,\right)\!\!\Bigg]\!.
$$
\end{theorem}

\section{Distribution related to $\tau_a$}

In this section, $N$ is assumed to be an even integer. The reader is referred
to~\cite{2006,2007,2008,ns,nish1,nish2,nish3}.

\subsection{Distribution of $(\tau_a,X(\tau_a))$}

Using the definition
$$
\mathbb{E}_x\!\!\left(e^{-\lambda\tau_a+i\mu X(\tau_a)}\right)
=\lim_{n\to+\infty}\mathbb{E}_x\!\!\left(e^{-\lambda \frac{k}{2^n}
+i\mu X(\frac{k}{2^n})}\ind_{\{X(\frac{k-1}{2^n})<a\le X(\frac{k}{2^n})\}}\right)\!,
$$
It can be seen that the Laplace-Fourier transform of the vector
$(\tau_a,X(\tau_a))$ is related to the pseudo-distribution of the vector
$(X(t),M(t))$ according as, for $\lambda>0$, $\mu\in\mathbb{R}$ and $x\le a$:
$$
\mathbb{E}_x\!\!\left(e^{-\lambda\tau_a+i\mu X(\tau_a)}\right)
=\left(\lambda+\mu^N\right) \int_0^{+\infty}e^{-\lambda t}
\,\mathbb{E}_x\!\!\left(e^{i\mu X(t)}\ind_{\{M(t)>a\}}\right)\mathrm{d}t.
$$
From this, it can be deduced that, for $x\le a$,
$$
\mathbb{E}_x\!\!\left(e^{-\lambda\tau_a+i\mu X(\tau_a)}\right)
=\sum_{j\in J} A_j\prod_{\ell\in J\setminus\{j\}}
(1-\frac{i\mu}{\!\sqrt[N\!]{\lambda}}\,\bar{\theta}_\ell)
\,e^{\theta_j\!\!\sqrt[N\!]{\lambda}\,(x-a)} \,e^{i\mu a}.
$$
In particular,
$$
\mathbb{E}_x\!\!\left(e^{-\lambda\tau_a}\right)
=\lambda \int_0^{+\infty} e^{-\lambda t}\,\mathbb{P}_{\!x} \{M(t)>a\}\,\mathrm{d}t
=\sum_{j\in J} A_j \,e^{\theta_j\!\!\sqrt[N\!]{\lambda}\,(x-a)}.
$$
The variables $\tau_a$ and $M(t)$ are related together according to
$\mathbb{P}_{\!x}\{\tau_a\le t\}=\mathbb{P}_{\!x}\{M(t)\ge a\}$.
The two-parameters Laplace-Fourier transform can be inverted (\cite{2007}).

\begin{theorem}[Lachal, 2006]
The joint pseudo-distribution of $(\tau_a,X(\tau_a))$ is given, for $x<a$, by
$$
\mathbb{P}_{\!x}\{\tau_a\in \mathrm{d}t,X(\tau_a)\in \mathrm{d}z\}/
\mathrm{d}t\,\mathrm{d}z=\sum_{p=0}^{N/2-1} \mathcal{J}_p(t;x-a)\,\delta_a^{(p)}(z)
$$
with $\mathcal{J}_p(t;\xi)={\displaystyle\mathop{\textstyle{\sum}}_{j\in J}}
\overline{\sigma_{j,p}} \,A_j \,I_p(t;\theta_j\xi)$. In particular,
$$
\mathbb{P}_{\!x}\{\tau_a\in \mathrm{d}t\}/\mathrm{d}t=\mathcal{J}_0(t;x-a),\quad
\mathbb{P}_{\!x}\{X(\tau_a)\in \mathrm{d}z\}/\mathrm{d}z
=\sum_{p=0}^{N/2-1} (-1)^p\frac{(x-a)^p}{p!}\,\delta_a^{(p)}(z).
$$
\end{theorem}
The $\delta_a^{(p)}$ are the successive derivatives of the Schwartz distribution
$\delta_a$, that is, for any test function $\phi$,
$<\delta_a^{(p)},\phi>=(-1)^p\phi^{(p)}(a)$.
The pseudo-distribution of $X(\tau_a)$ is remarkable since it means that the
pseudo-process $(X(t))_{t\ge 0}$ is formally concentrated at the site $a$
at time $\tau_a$ in a ``distributional'' sense.

\begin{example}
In the case $N=4$, Nishioka (\cite{nish1,nish2}) obtained the remarkable result
$$
\mathbb{P}_{\!x}\{X(\tau_a)\in \mathrm{d}z\}/\mathrm{d}z=\delta_a(z)-(x-a)\delta_a'(z).
$$
Moreover,
$$
\mathbb{P}_{\!x}\{\tau_a\in \mathrm{d}t,X(\tau_a)\in \mathrm{d}z\}/
\mathrm{d}t\,\mathrm{d}z
=\mathcal{J}_0(t;x-a)\,\delta_a(z)+\mathcal{J}_1(t;x-a)\,\delta'_a(z)
$$
with
\begin{align*}
\mathcal{J}_0(t;\xi)
&=\frac{\xi}{2\pi t} \int_0^{+\infty} \left(e^{\xi\lambda}
-\cos(\xi\lambda)+\sin(\xi\lambda)\right) e^{-t\lambda^4}\,\mathrm{d}\lambda,
\\
\mathcal{J}_1(t;\xi)
&=\frac{2}{\pi} \int_0^{+\infty} \left(\cos(\xi\lambda)+\sin(\xi\lambda)
-e^{\xi\lambda}\right) \lambda^2\,e^{-t\lambda^4}\,\mathrm{d}\lambda.
\end{align*}
\end{example}

\subsection{Distribution of $(\tau_a^b,X(\tau_a^b))$}

The pseudo-distributions of the vectors $(\tau_a^b,X^b(\tau_a^b))$ and
$(X^b(t),M^b(t))$ are related together according as, for $\lambda>0$,
$\mu\in\mathbb{R}$ and $x\le a$,
$$
\mathbb{E}_x\!\!\left(e^{-\lambda\tau_a^b+i\mu X^b(\tau_a^b)}
\ind_{\{\tau_a^b<\infty\}}\right)=\left(\lambda-ib\mu+\mu^N\right)
\int_0^{+\infty}e^{-\lambda t}\,\mathbb{E}_x\!\!\left(e^{i\mu X^b(t)}
\ind_{\{M^b(t)>a\}}\right)\mathrm{d}t
$$
from which it comes that the Laplace-Fourier transform of the vector $(\tau_a^b,X^b(\tau_a^b))$
is given, for $\lambda>0$, $\mu\in\mathbb{R}$ and $x\le a$, by
$$
\mathbb{E}_x\!\!\left(e^{-\lambda\tau_a^b+i\mu X^b(\tau_a^b)}
\ind_{\{\tau_a^b<\infty\}}\right)=\sum_{j\in J} \left(\prod_{
\ell\in J\setminus\{j\}}\frac{\omega_\ell^b(\lambda)+\mu}
{\omega_\ell^b(\lambda)-\omega_j^b(\lambda)}\right)
e^{i\mu a-i\omega_j^b(\lambda)(x-a)}.
$$
\begin{example}
In the case $N=4$, Nakajima and Sato (\cite{ns}) derived the following result:
$$
\mathbb{E}_x\!\!\left(e^{-\lambda\tau_a^b+i\mu X^b(\tau_a^b)}
\ind_{\{\tau_a^b<\infty\}}\right)=e^{i\mu a}
\left(\frac{\mu+\omega_1}{\omega_2-\omega_1}\,e^{-i\omega_1(x-a)}
+\frac{\mu+\omega_2}{\omega_1-\omega_2}\,e^{-i\omega_2(x-a)}\right)
$$
where $\omega_1$ and $\omega_2$ are the two roots of the polynomial
$X^4+ibX+\lambda$ having positive imaginary part.
\end{example}

We can deduce the pseudo-distribution of the overshooting place
$X^b(\tau_a^b)$ on the set $\{\tau_a^b<\infty\}$ and the pseudo-probability
of eventually hitting the interval $[a,+\infty)$ (\cite{2008}).

\begin{theorem}[Lachal, 2006]
The pseudo-distribution of the overshooting place
$X^b(\tau_a^b)\times\ind_{\{\tau_a^b<\infty\}}$ is given, for $x\le a$, by
$$
\mathbb{P}_{\!x}\{X^b(\tau_a^b)\in \mathrm{d}z,\tau_a^b<\infty\}/\mathrm{d}z=
\sum_{p=0}^{N/2-1} \frac{i^p}{|b|^{\frac{p}{N-1}}} \Bigg(\sum_{j\in \tilde{J}}
\frac{\tilde{\sigma}_{j,N/2-p-1}}{\displaystyle
\mathop{\textstyle{\prod}}_{\ell\in\tilde{J}\setminus\{j\}}
(\tilde{\theta}_\ell-\tilde{\theta}_j)} e^{-i\tilde{\theta}_j|b|^{\frac{1}{N-1}}
(x-a)}\Bigg) \delta_a^{(p)}(z)
$$
where the $\tilde{\sigma}_{j,p}$'s are given by $\tilde{\sigma}_{j,0}=1$
and for $1\le p\le N/2-1$,
$$
\tilde{\sigma}_{j,p}=\sum_{\ell_1<\cdots<\ell_p\atop \ell_1,\ldots,
\ell_p\in \tilde{J}\setminus\{j\}}\tilde{\theta}_{\ell_1}\ldots \tilde{\theta}_{\ell_p}.
$$
The pseudo-probability of eventually overshooting the level $a$ is given by
$$
\mathbb{P}_{\!x}\{\tau_a^b<\infty\}=
\begin{cases}
1 & \mbox{if $b>0$,}
\\
\displaystyle\sum_{j\in \tilde{J}}
\frac{e^{-i\tilde{\theta}_j|b|^{\frac{1}{N-1}}(x-a)}}{\displaystyle
\mathop{\textstyle{\prod}}_{\ell\in \tilde{J}\setminus\{j\}}
\Big(1-e^{i\frac{2(j-\ell)}{N-1}\pi}\Big)}& \mbox{if $b<0$.}
\end{cases}
$$
\end{theorem}
\begin{example}
For $N=4$, the above results yield for $x\le a$, in the case where $b>0$,
$$
\mathbb{P}_{\!x} \{\tau_a^b<\infty\}=1,\quad
\mathbb{P}_{\!x} \{X^b(\tau_a^b)\in \mathrm{d}z,\tau_a^b<\infty\}/\mathrm{d}z
=\delta_a(z)+\frac{1}{\sqrt[3]b}\big(1-e^{\!\sqrt[3]b\,(x-a)}\big)\delta_a'(z),
$$
and, in the case where $b<0$,
$$
\mathbb{P}_{\!x} \{\tau_a^b<\infty\}=p_0^b(x-a),\quad
\mathbb{P}_{\!x} \{X^b(\tau_a^b)\in \mathrm{d}z,\tau_a^b<\infty\}/\mathrm{d}z
=p_0^b(x-a)\delta_a(z)+p_1^b(x-a)\delta_a'(z)
$$
with
$$
p_0^b(\xi)=\frac{2}{\sqrt3}\,e^{\frac12\sqrt[3]{|b|}\,\xi}
\cos\!\bigg(\frac{\sqrt3}{2}\,\sqrt[3]{|b|}\,\xi+\frac{\pi}{6}\bigg)\!,
\quad p_1^b(\xi) =-\frac{2}{\sqrt3\,\sqrt[3]{|b|}}\,e^{\frac12\sqrt[3]{|b|}\,\xi}
\sin\!\bigg(\frac{\sqrt3}{2}\,\sqrt[3]{|b|}\,\xi\bigg)\!.
$$
\end{example}

\section{Works in progress}

The problem of the two-sided barrier $\{a,b\}$ is much more difficult to tackle
than the single one. The variables $T_{ab}(t)$, $\tau_{ab}$ and the
maximum/minimum functionals are related together according as, for $x\in(a,b)$,
$$
\mathbb{P}_{\!x}\{\tau_{ab}\ge t\}=\mathbb{P}_{\!x}\{a\le m(t)\le M(t)\le b\}
=\mathbb{P}_{\!x}\{T_{ab}(t)=t\}
=\lim_{\mu\to +\infty} \mathbb{E}_x\!\!\left(e^{-\mu (T_{ab}(t)-t)}\right)\!.
$$

\subsection{Distributions related to $T_{ab}(t)$}

The variable $T_{ab}(t)$ is introduced in \cite{bo} for defining a local time
for the pseudo-process $(X(t))_{t\ge0}$. Set
$\Phi(x)=\int_0^{+\infty} e^{-\lambda t} \, \mathbb{E}_x\!\!
\left(e^{-\mu T_{ab}(t)}\right) \mathrm{d}t$ for $\lambda,\mu>0$ and
$x\in\mathbb{R}$. The function $\Phi$ satisfies the system
$$
\kappa_{_{\! N}}\frac{\mathrm{d}^N \Phi}{\mathrm{d} x^N}(x)
=\left\{\begin{array}{ll}
(\lambda+\mu)\,\Phi(x)-1  & \mbox{for } x\in(a,b),
\\[1ex]
\lambda\, \Phi(x)-1 & \mbox{for } x\notin(a,b),
\end{array}\right.
$$
and
$$
\forall k\in\{0,1,\ldots,N-1\},\;\frac{\mathrm{d}^k \Phi}{\mathrm{d} x^k}(a^+)
=\frac{\mathrm{d}^k \Phi}{\mathrm{d} x^k}(a^-)\mbox{ and }
\displaystyle\frac{\mathrm{d}^k \Phi}{\mathrm{d} x^k}(b^+)
=\frac{\mathrm{d}^k \Phi}{\mathrm{d} x^k}(b^-).
$$
It seems to be difficult to solve explicitly this system.
In the same way, the problem of the joint pseudo-distribution of
$(X(t),T_{ab}(t))$ should be more complicated.

The sojourn time within a strip is used by Beghin \& Orsingher (\cite{bo})
for defining a local time at 0 for $(X(t))_{t\ge 0}$ :
$L(t)=\lim_{\varepsilon\to 0^+}\frac{1}{\varepsilon}
\int_0^t\ind_{\{X(s)\in[-\varepsilon,\varepsilon]\}}\, \mathrm{d}s$.
They obtained the very simple formula
$$
\int_0^{+\infty} e^{-\lambda t} \, \mathbb{E}_0\!(e^{-\mu L(t)}) \,\mathrm{d}t
= \frac{1}{\lambda+c\mu\!\sqrt[N\!]{\lambda}}
$$
where
$$
c=\frac{1}{\alpha N\sin\frac{\pi}{\alpha N}},\;\textstyle \alpha=
\left\{\begin{array}{l}
1 \mbox{ if $N$ is odd,}\\
2 \mbox{ if $N$ is even.}
\end{array}\right.
$$

\subsection{Distributions related to $\tau_{ab}$}

By introducing a pseudo-random walk as in \cite{sato} and studying the
similar functional to $\tau_{ab}$ (first exit time from a finite interval),
we can derive the pseudo-distribution of $(\tau_{ab},X(\tau_{ab}))$
in the case where $N$ is even. We have obtained the following result
(\cite{2010}).
\begin{theorem}[Lachal, 2010]
The pseudo-distribution of $X(\tau_{ab})$ has the form
$$
\mathbb{P}_{\!x}\{\tau_{ab}\in \mathrm{d}t,X(\tau_{ab})\in \mathrm{d}z\}/
\mathrm{d}t\,\mathrm{d}z=\sum_{p=0}^{N/2-1} \mathcal{J}_p(x)\,\delta_a^{(p)}(z)
+\sum_{p=0}^{N/2-1} \mathcal{K}_p(x)\,\delta_b^{(p)}(z)
$$
where $\mathcal{J}_p$ and $\mathcal{K}_p$ are some functions.
The pseudo-distribution of $X(\tau_{ab})$ is given by
$$
\mathbb{P}_{\!x}\{X(\tau_{ab})\in \mathrm{d}z\}/\mathrm{d}z
=\sum_{p=0}^{N/2-1} H_p^-(x)\,\delta_a^{(p)}(z)+\sum_{p=0}^{N/2-1} H_p^+(x)
\,\delta_b^{(p)}(z)
$$
where the functions $H_p^-$ and $H_p^+$, $0\le p\le N/2-1$, are the interpolation
Hermite polynomials such that
$$
\frac{\mathrm{d}^qH^-_p}{\mathrm{d}x^q}(a)=\delta_{pq},
\;\frac{\mathrm{d}^qH^-_p}{\mathrm{d}x^q}(b)=0,\;
\frac{\mathrm{d}^qH^+_p}{\mathrm{d}x^q}(a)=0,
\;\frac{\mathrm{d}^qH^+_p}{\mathrm{d}x^q}(b)
=\delta_{pq} \mbox{ for $0\le q\le N/2-1$.}
$$
In particular, the ``ruin pseudo-probabilities'' are given by
$$
\mathbb{P}_{\!x}\{\tau_a^-<\tau_b^+\}=H_0^-(x)
\mbox{ and }\mathbb{P}_{\!x}\{\tau_b^+<\tau_a^-\}=H_0^+(x).
$$
\end{theorem}
\begin{example}
In the case $N=4$, the above results supply
$$
\mathbb{P}_{\!x}\{X(\tau_{ab})\in \mathrm{d}z\}/\mathrm{d}z
=H_0^-(x)\,\delta_a(z)+H_1^-(x)\,\delta_a'(z)+H_0^+(x)\,\delta_b(z)
+H_1^+(x)\,\delta_b'(z)
$$
where
$$
H_0^-(x)=\frac{(x-b)^2(2x+b-3a)}{(b-a)^3},\quad
H_1^-(x)=-\frac{(x-a)(x-b)^2}{(b-a)^2},
$$
\vspace{-2ex}
$$
H_0^+(x)=-\frac{(x-a)^2(2x+a-3b)}{(b-a)^3},\quad
H_1^+(x)=-\frac{(x-a)^2(x-b)}{(b-a)^2}.
$$
\end{example}

\subsection{Distributions related to $\sigma_{ab}(t)$}

By applying the pseudo-Markov property, it can be easily seen that the
distributions of $\sigma_{ab}(t)$ and $\tau_{ab}$ are related together
according as, for $x\in\mathbb{R}$ and $\sigma\in[0,t]$,
$$
\mathbb{P}_{\!x}\{\sigma_{ab}(t)\le\sigma\}=\mathbb{P}_{\!x}\{\forall
s\in[\sigma,t], X(s)\in [a,b]\}=\mathbb{E}_x\big(\mathbb{P}_{X(\sigma)}
\{\tau_{ab}\ge t-\sigma\}\ind_{\{X(\sigma)\in[a,b]\}}\big).
$$
This equality may be extended to an ``excursion'' between the thresholds $a$
and $b$. Indeed, by introducing $\varsigma_{ab}(t)=\inf\{s\ge t:X(t)\notin (a,b)\}$,
we have, for $x\in\mathbb{R}$ and $0\le\sigma\le t \le\varsigma$,
\begin{align*}
\mathbb{P}_{\!x}\{\sigma_{ab}(t)\le\sigma,\varsigma_{ab}(t)\ge\varsigma,
X(\varsigma_{ab}(t))\in \mathrm{d}y\}
&
=\mathbb{P}_{\!x}\{\forall s\in[\sigma,\varsigma], X(s)\in [a,b],
X(\varsigma_{ab}(t))\in \mathrm{d}y\}
\\
&
=\mathbb{E}_x\big(\mathbb{P}_{X(\sigma)}\{\tau_{ab}\ge\varsigma-\sigma,
X(\tau_{ab})\in \mathrm{d}y\}\ind_{\{X(\sigma)\in[a,b]\}}\big).
\\
&
=\int_a^b p(\sigma;x-z)\,\mathbb{P}_{\!z}\{\tau_{ab}\ge\varsigma-\sigma,
X(\tau_{ab})\in \mathrm{d}y\}\big)\,\mathrm{d}z.
\end{align*}
Acknowledgements.
This survey is the text associated with the talk I gave at the fifth IWAP held
in Madrid in July 2010 at the session ``special stochastic processes'' chaired
by K.J.~Hochberg. It has been a great honor for me to give this talk in the
presence of K.J.~Hochberg and E.~Orsingher who were pioneers in the domain of
pseudo-processes.


\end{document}